\documentclass[11pt,leqno]{amsart}
\usepackage{amsthm,amsfonts,amssymb,amsmath,oldgerm}
\numberwithin{equation}{section}
\usepackage{graphicx}
\usepackage{fullpage}

\renewcommand\d{\partial}
\newcommand\dd{{\textrm d}}

\renewcommand\d{\partial}

\newcommand\R{\mathbb R}
\newcommand\C{\mathbb C}

\newcommand\br{\begin{remark}}
\newcommand\er{\end{remark}}
\newcommand\bp{\begin{pmatrix}}
\newcommand\ep{\end{pmatrix}}
\newcommand{\be}{\begin{equation}}
\newcommand{\ee}{\end{equation}}
\newcommand\ba{\begin{equation}\begin{aligned}}
\newcommand\ea{\end{aligned}\end{equation}}

\newcommand{\bap}{\begin{app}}
\newcommand{\eap}{\end{app}}
\newcommand{\begs}{\begin{exams}}
\newcommand{\eegs}{\end{exams}}
\newcommand{\beg}{\begin{example}}
\newcommand{\eeg}{\end{exaplem}}
\newcommand{\bpr}{\begin{proposition}}
\newcommand{\epr}{\end{proposition}}
\newcommand{\bt}{\begin{theorem}}
\newcommand{\et}{\end{theorem}}
\newcommand{\bc}{\begin{corollary}}
\newcommand{\ec}{\end{corollary}}
\newcommand{\bl}{\begin{lemma}}
\newcommand{\el}{\end{lemma}}
\newcommand{\bd}{\begin{definition}}
\newcommand{\ed}{\end{definition}}
\newcommand{\brs}{\begin{remarks}}
\newcommand{\ers}{\end{remarks}}

\newcommand{\const}{\text{\rm constant}}
\newcommand{\Id}{{\rm Id }}

\newtheorem{theorem}{Theorem}[section]
\newtheorem{proposition}[theorem]{Proposition}
\newtheorem{corollary}[theorem]{Corollary}
\newtheorem{lemma}[theorem]{Lemma}

\theoremstyle{remark}
\newtheorem{remark}[theorem]{Remark}
\theoremstyle{definition}
\newtheorem{definition}[theorem]{Definition}

\newtheorem{example}[theorem]{Example}

\newcommand{\RM}{\mathbb{R}}

\newcommand{\f}{\frac}

\newcommand{\beq}{\begin{equation}}
\newcommand{\eeq}{\end{equation}}

\title{
	Periodic-coefficient damping estimates, and 
	stability of large-amplitude roll waves in inclined thin film flow
}

\author{L.Miguel Rodrigues}
\address{Universit\'e Lyon 1, Institut Camille Jordan, INRIA \'EP Kaliffe, Villeurbanne, France}
\email{{\tt rodrigues@math.univ-lyon1.fr}}
\thanks{Research of L.M.R. was partially supported by the ANR project BoND ANR-13-BS01-0009-01.}
\author{Kevin Zumbrun}
\address{Indiana University, Bloomington, IN 47405}
\email{{\tt kzumbrun@indiana.edu}}
\thanks{Research of K.Z. was partially supported under NSF grant no. DMS-0300487.}
\begin{document}

\begin{abstract}
A technical obstruction preventing the conclusion of nonlinear stability of large-Froude number roll waves of the St. Venant
equations for inclined thin film flow is the "slope condition" of Johnson-Noble-Zumbrun, used to obtain pointwise symmetrizability of the linearized equations and thereby high-frequency resolvent bounds and a crucial $H^s$ nonlinear damping estimate. Numerically, this condition is seen to hold for Froude numbers $2<F\lessapprox 3.5$, but to fail for $3.5 \lessapprox F$.  As hydraulic engineering applications typically involve Froude number $3\lessapprox F \lessapprox 5$, this issue is indeed relevant to practical considerations. Here, we show that the pointwise slope condition can be replaced by an averaged version which holds always, thereby completing the nonlinear theory in the large-$F$ case. The analysis has potentially larger interest as an extension to the periodic case of a type of weighted ``Kawashima-type'' damping estimate introduced in the asymptotically-constant coefficient case for the study of stability of large-amplitude viscous shock waves.
\end{abstract}
\date{\today}
\maketitle

\section{Introduction}\label{s:introduction}

The St. Venant equations of inclined thin film flow, in nondimensional Lagrangian form, are
	\ba
	\label{swl}
	\partial_t \tau-\partial_x u&=0,
	\\
	\partial_t u+\partial_x\left(\frac{\tau^{-2}}{2F^2}\right)&=1-\tau\,u^2+\nu\partial_x(\tau^{-2}\partial_x u),
\ea
where $\tau=1/h$ is the reciprocal of fluid height $h$, $u$ is tangential fluid velocity averaged with respect to height,
$x$ is a Lagrangian marker,  
$F$ is a Froude number given by the ratio between a chosen reference speed of the fluid 
and speed of gravity waves, and $\nu=R_e^{-1}$, with $R_e$ the Reynolds number of the fluid.
The terms $1$ and $\tau u^2$ on the righthand side of the second
equation model, respectively, gravitational force and turbulent friction along the bottom.
Roughly speaking, $F$ measures inclination, with $F=0$ corresponding to horizontal
and $F\to \infty$ to vertical inclination of the plane.

An interesting and much-studied phenomenon in thin film flow is the appearance of
{\it roll-waves}, or spatially periodic traveling-waves corresponding to solutions
\be\label{roll}
(\tau, u)(x,t)=(\bar \tau, \bar u)(x-ct)
\ee
of \eqref{swl}. 
These are well-known {\it hydrodynamic instabilities,} 
arising for \eqref{swl} in the region $F>2$ for which constant solutions, 
corresponding to parallel flow, are unstable,
with applications to landslides, river and spillway flow, and topography of
sand dunes and sea beds \cite{BM}. 

Nonlinear stability of roll-waves themselves has been a long-standing open problem.
However, this problem has recently been mostly solved in a series of works
by the authors together with Barker, Johnson, and Noble; see
\cite{JZN,BJRZ,BJNRZ1,JNRZ,BJNRZ2}.
More precisely, it has been shown that, under a certain technical condition having to do 
with the slope of the traveling-wave profile $(\bar \tau, \bar u)$,
{\it spectral stability} in the sense of Schneider \cite{S1,S2,JZN,JNRZ},
{\it implies linear and nonlinear modulational stability with optimal rates of decay}, and,
moreover, asymptotic behavior is well-described by a system of second-order {\it Whitham equations}
obtained by formal WKB expansion.

In turn, spectral stability has been characterized analytically in the {\it weakly unstable} limit $F\to 2$
and numerically for intermediate to large $F$ in terms of two simple power-law descriptions, 
in the small- and large-$F$ regimes, respectively,
of the band of periods $X$ for which roll waves are spectrally stable, 
as functions of $F$ and discharge rate $q$ (an invariant of the flow
describing the flux of fluid through a given reference point) \cite{BJNRZ2}.
That is, apart from the technical slope condition, there is at this point a rather complete theory of spectral, linear,
and nonlinear stability of roll wave solutions of the St. Venant equations.
However, up to now it was not clear whether failure of the slope condition was a purely
technical issue or might be an additional mechanism for instability.

Precisely, this slope condition reads, in the Lagrangian formulation \eqref{swl}--\eqref{roll}, as
\be\label{e:Lslope}
2 \nu \bar u_x < F^{-2},
\ee
where $\bar u$ is the velocity component of traveling wave \eqref{roll}.
It is seen numerically to be satisfied for $F\lessapprox 3.5$, but to fail for $F\gtrapprox 3.5$ \cite{BJNRZ2}.
For comparison, hydraulic engineering applications typically involve Froude numbers $2.5\lessapprox F\lessapprox 20$
\cite{A,Br1,Br2}; {\it hence \eqref{e:Lslope} is a real 
physical 
restriction.} 
From the mathematical point of view, the distinction is between small-amplitude, slowly varying waves for which \eqref{e:Lslope}
is evidently satisfied and large-amplitude, rapidly-varying waves, such as appear in the spectrally stable regime for 
small and large $F$, respectively \cite{BJNRZ2}.

\medskip

The 
role of condition \eqref{e:Lslope} in the stability analysis is
to obtain pointwise symmetrizability of the linearized equations 
and thereby high-frequency resolvent bounds and a crucial 
nonlinear damping estimate used to control higher derivatives in a nonlinear iteration scheme.
The purpose of the present brief note is to show, by a refined version of the energy estimates of
\cite{JZN,BJRZ}, that the pointwise condition \eqref{e:Lslope} can be replaced by
an {\it averaged version} that is {\it always satisfied,} 
while still retaining the high-frequency resolvent and nonlinear damping estimates needed for the nonlinear analysis of
\cite{JZN,JNRZ},
thus {\it effectively completing the nonlinear stability theory.}

\medskip

The remainder of this paper is devoted to establishing the requisite weighted energy estimates,
first, in Sections \ref{s:prelim}-\ref{s:damp},
in the simplest, linear time-evolution setting 
then, in Sections \ref{s:hf} and \ref{s:ndamp}, respectively,
in the closely related high-frequency resolvent and nonlinear time-evolution settings.
The estimates so derived may be seen to be periodic-coefficient analogs of weighted ``Kawashima-type'' estimates
derived in the asymptotically-constant coefficient case
for the study of stability of large-amplitude viscous shock waves \cite{Z1,Z2,GMWZ},
to our knowledge the first examples of such estimates specialized to the periodic setting.
We discuss this connection in Sections \ref{s:shock} \& \ref{s:discussion}. 
More, this
seems to be the first instance of a proof of hypocoercive\footnote{The reader interested in replacing Kawashima-type estimates in the more general context of hypocoercive decay estimates is referred to \cite[Remark~17]{Villani} and references in \cite[Appendix~A]{R}, especially \cite{Beauchard-Zuazua}.} decay where periodicity is used in a crucial way. 
We note, finally, the relation between these weights and the
``gauge functions'' used for similar purposes in short-time (i.e., well-posedness) dispersive theory \cite{LP,BDD,Mietka},
a connection brought out further by our choice of notation in the proof.
This indicates perhaps a potential for wider applications of these ideas in the study of periodic wave trains.

\section{Preliminary observations}\label{s:prelim}
Making the change of variables $x\to x-ct$ to
co-moving coordinates, we convert \eqref{swl} to
\ba \label{eqn:co1conslaw}
\tau_t-c\tau_x - u_x&= 0,\\
u_t-cu_x + ((2F^2)^{-1}\tau^{-2})_x&=
1- \tau u^2 +\nu (\tau^{-2}u_x)_x ,
\ea
and the traveling-wave solution
to a stationary solution $U(x,t)=(\tau(x,t),u(x,t))=(\bar \tau(x),\bar u(x))$
convenient for stability analyis.

We note for later that the traveling-wave ODE becomes
\be \label{ODE}
-c\bar \tau_x - \bar u_x= 0,\qquad
-c\bar u_x + ((2F^2)^{-1}\bar \tau^{-2})_x=
1- \bar \tau \bar u^2 +\nu (\bar \tau^{-2}\bar u_x)_x ,
\ee
yielding the key fact that
\be\label{fact}
f(\bar \tau)\bar u_x= cf(\bar \tau)\bar \tau_x
\ee
is a perfect derivative for any function $f(\cdot)$, hence {\it zero mean} over one period.
We note also as in \cite{JZN} that $c\neq0$, else $u\equiv \const$ and the equation for $\tau$
reduces to first order, hence does not admit nontrivial periodic solutions.
Linearizing about $\bar U=(\bar \tau,\bar u)$ gives the {\it linearized equations}
\ba \label{e:lin}
\tau_t-c\tau_x - u_x&= 0,\\
u_t-cu_x -(\alpha \tau)_x 
&=
\nu (\bar \tau^{-2}u_x)_x 
- \bar u^2 \tau 
- 2\bar u \bar \tau u ,
\ea
where
\be\label{alpha}
\alpha:= \bar \tau^{-3}(F^{-2} + 2\nu\bar u_x ).
\ee

With this notation, the slope condition of \cite{JZN} appears as $\bar \tau^3 \alpha>0$.
We note that, by \eqref{fact}, the mean over one period of $g(\bar \tau)\alpha$ is {\it positive} for any positive $g$:
\be\label{pos}
\langle g(\bar \tau) \alpha \rangle= \langle g(\bar \tau)\bar \tau^{-3} F^{-2} \rangle >0.
\ee
That is, (any reasonable version of) the slope condition holds always in an averaged sense.\footnote{
Here and elsewhere we use $\langle h \rangle$ to denote mean over one period of a function $h$.}
An approximate asymptotic diagonalization in the large spectrum regime-- see \cite{BJRZ,BJNRZ2}, 
in particular \cite[Appendix~A]{BJNRZ2}-- reveals that the sharp\footnote{In the sense that there exist curves of spectrum for the operator $L$-- defined below-- that are going to infinity and whose real parts converge to $-\langle \alpha \bar \tau^2\rangle/\nu$.} relevant averaged conditions is 
$$
\frac{\langle \alpha \bar \tau^2\rangle}{\nu} >0\,.
$$
We shall show in the rest of the paper that this averaged condition is in fact sufficient for the nonlinear analysis of \cite{JZN,JNRZ}.

\section{Linear damping estimate }\label{s:damp}
Introduce now some `gauge' functions $\phi_1$, $\phi_2$ and $\phi_3$ and define for $U=(\tau,u)$ the energy
\be\label{form}
\mathcal{E}(U):=
\int \Big(
\tfrac12\phi_1 \tau_x^2 +
\tfrac12\phi_2\bar \tau^3  u_x^2 +
\phi_3 \tau u_x \Big).
\ee

A brief computation yields that solutions $U$ of \eqref{e:lin} satisfy
\ba\label{diss}
\frac{d}{dt}\mathcal{E}(U(t))&=
\int\Big(
-\left(\tfrac c 2 (\phi_1)_x + \alpha \phi_3\right)\tau_x^2
-\left(\tfrac{\nu} {\bar \tau^2} \phi_2\right)u_{xx}^2
+ \left(\phi_1 -\alpha\phi_2 + \tfrac{\nu}{\bar \tau^2} \phi_3\right) \tau_x u_{xx} \Big)\\
&\quad
+ O\big((\|u\|_{H^2}+\|\tau\|_{H^1}) (\|u\|_{H^1}+\|\tau\|_{L^2})\big).
\ea
The original gaugeless strategy that works when $\alpha$ is positive may be achieved by choosing $\phi_1\equiv 1$, $\phi_2=\phi_1/\alpha$ and $0<\phi_3\equiv \const \ll 1$. The possibility of choosing $\phi_2=\phi_1/\alpha$ while keeping both $\phi_1$ and $\phi_2$ positive is a direct manifestation of the fact that in this case the first-order part of system~\eqref{e:lin} is symmetrizable. For the general case, of interest here, we instead take
\ba\label{take1}
\frac c 2 (\phi_1)_x + \Big(\frac{\alpha \bar \tau^2}{\nu}
-\frac{\langle \alpha \bar \tau^2\rangle}{\nu}\Big)\phi_1&=0,\qquad\phi_1(0)>0,\\
\ea
\ba\label{take2}
\phi_1 -\alpha\phi_2 + \frac{\nu}{\bar \tau^2} \phi_3&=0,\qquad 0<\phi_2\equiv \const \ll 1,\\
\ea
so that $\phi_3$ is chosen to kill the 
indefinite
cross-term and the fact that $\phi_1$ is not constant and thus does not commute with the generator of system~\eqref{e:lin} is used to average and cancel the 
``bad''
oscillating part of $\frac{\alpha \bar \tau^2}{\nu}$ through the arising nontrivial commutator.

With these choices,
we obtain after another brief computation
\ba\label{gooddiss}
\frac{d}{dt}\mathcal{E}(U(t))&=
-\int
\Big[ \left(\tfrac{\langle\bar \tau^2 \alpha\rangle}{\nu} \phi_1 -\tfrac{\alpha^2\bar \tau^2}{\nu}\phi_2\right)\tau_x^2
+\left(\tfrac{\nu}{\bar \tau^2} \phi_2\right)u_{xx}^2
\Big]\\
&\quad
+ O\big(\|u\|_{H^2}+\|\tau\|_{H^1}) (\|u\|_{H^1}+\|\tau\|_{L^2})\big)\\
&\leq
-\eta_1 (\|u_{xx}\|_{L^2}^2 + \|\tau_x\|_{L^2}^2)
+ C_1 \big(\|u\|_{H^2}+\|\tau\|_{H^1}) (\|u\|_{H^1}+\|\tau\|_{L^2})\big),
\ea
for some positive $\eta_1$ and $C_1$, whence, by interpolation inequality $\|u\|_{H^1}\lesssim \|u\|_{L^2}^{1/2}\|u\|_{H^2}^{1/2}$ and the fact that $\mathcal{E}(U)\sim 
(\|u_{x}\|_{L^2}^2 + \|\tau_x\|_{L^2}^2)$ modulo $\|\tau\|_{L^2}^2$,
\be\label{dampest}
\frac{d}{dt}\mathcal{E}(U(t))\leq -\eta \mathcal{E}(U(t)) + C\|U(t)\|_{L^2}^2,
\ee
for some positive $\eta$ and $C$, a standard {\it linear damping estimate}.

Note that, in the step $\mathcal{E}(U)\sim 
(\|u_{x}\|_{L^2}^2 + \|\tau_x\|_{L^2}^2)$ modulo $L^2$,
we have used in a critical way that 
$$\int^x\Big(\frac{\alpha \bar \tau^2}{\nu}-\frac{\langle \alpha \bar \tau^2\rangle}{\nu}\Big),$$
hence $\phi_1$ and $1/\phi_1$, remains bounded, a consequence of periodicity plus zero mean.

To derive a first corollary from the key estimate \eqref{dampest}, one may combine it with the standard $L^2$ bound 
$$
\frac{d}{dt}\int \Big(\tfrac12\tau^2 +\tfrac12u^2 \Big)\ =\ 
-\int\tfrac{\nu} {\bar \tau^2}u_{x}^2
\ +\ O\big((\|u\|_{H^1}+\|\tau\|_{L^2}) (\|u\|_{L^2}+\|\tau\|_{L^2})\big)
$$
to obtain the following lemma.

\bl\label{l:sample}
There exist positive $\theta$ and $C$ such that any $U$ solving \eqref{e:lin} satisfies for any $t\geq0$
$$
\|U(t)\|_{H^1}\ \leq\ C\,e^{-\theta\,t} \|U(0)\|_{H^1}\ +\ C\,\int_0^t e^{-\theta\,(t-s)} \|U(s)\|_{L^2} \dd s\,.
$$
\el

\section{Applications}\label{s:applications}

Lemma~\ref{l:sample} 
by itself
is not of much direct practical use. 
However, as we will now show,
we can readily adapt
its proof, and especially estimate~\eqref{dampest}, to obtain various useful forms of high-frequency damping estimates. The reader unfamiliar with
these considerations may benefit from first having a look at \cite[Appendix~A]{R} for a terse introduction to 
this approach.
Indeed what follows stems directly from the mere introduction 
in the classical strategy described there of gauges leading to \eqref{dampest}. 
See also \cite{Z1,Z2} for related estimates in the shock wave case.

\subsection{High-frequency resolvent bound}\label{s:hf}

An important part of the 
proofs in \cite{JZN,JNRZ} is dedicated to estimates of semigroups generated by linearization around a given wave, to be used in 
an
integral formulation of the original nonlinear systems. 
These estimates
are deduced from spectral considerations and the noncritical part of the linearized evolution is directly controlled by an abstract spectral gap argument that only requires 
uniform bounds on 
certain
resolvents. Our claim is that a spectral version of Lemma~\ref{l:sample} does provide these uniform bounds.  

To be more specific let $L$ denote the operator generating the linearized evolution around $\bar U$, that is, such that system \eqref{e:lin} reads $U_t-LU$. The operator $L$ is a differential operator with periodic coefficients but acting on functions defined on the full line. We do not apply directly spectral considerations to $L$ but rather to its operator-valued Bloch symbols $L_\xi$, associated with the Floquet-Bloch transform-- see \cite{JNRZ,R} for instance. Explicitly, if $\Xi$ denotes the fundamental period of $\bar U$, for any Floquet exponent $\xi$ in the Brillouin zone $[-\pi/\Xi, \pi/\Xi)$, the operator $L_\xi$ acts on functions of period $\Xi$ by $L_\xi:= e^{-i\xi \cdot}L e^{i\xi \cdot}$.

The operator $\d_x$ itself has Bloch symbols $\d_x+i\xi$. As a result, when dealing with $L_\xi$, the (equivalent) norm of interest on $H^s_{per}(0,\Xi)$ is $\|f\|_{H^s_\xi}=\left(\sum_{k=0}^s\|(\d_x+i\xi)^kf\|_{L^2(0,\Xi)}^2\right)^{1/2}$. Consider now the resolvent equation
\be\label{res}
(\lambda-L_\xi)U\ =\ F\,.
\ee
Letting $\langle\cdot,\cdot\rangle$ denote complex inner product, we find by computations
essentially identical to those in Section \ref{s:damp}, substituting $\lambda U$ for $U_t$ and $\d_x+i\xi$ for $\d_x$, that,
defining $\phi_j$ as in \eqref{take1}--\eqref{take2}, and
$$
\mathcal {E}_\xi(U):=
\tfrac12\langle \phi_1 (\d_x+i\xi)\tau,(\d_x+i\xi)\tau\rangle 
+\tfrac12\langle \phi_2 (\d_x+i\xi)u, (\d_x+i\xi)u \rangle 
+\Re \langle \phi_3 \tau, (\d_x+i\xi)u\rangle,
$$
one derives
\ba\label{re}
2\Re(\lambda)\ \mathcal{E}_\xi(U)& =\Re(2\lambda \mathcal{E}_\xi(U))\\
&=\ \Re (\langle \phi_1 (\d_x+i\xi)u,(\d_x+i\xi)(\lambda u)\rangle 
+\langle \phi_2 (\d_x+i\xi)\tau, (\d_x+i\xi)(\lambda \tau) \rangle\\[0.25em]
&\quad+ \langle \phi_3 \tau, (\d_x+i\xi)(\lambda u)\rangle
+\langle \phi_3 (\d_x+i\xi)u, \lambda \tau\rangle)\\
&=
-\int(\tfrac c 2 (\phi_1)_x + \alpha \phi_3)|(\d_x+i\xi)\tau|^2
\ -\int(\tfrac{\nu}{\bar \tau^2} \phi_2)|(\d_x+i\xi)^2u|^2\\
&\quad+\int(\phi_1 -\alpha\phi_2 + \tfrac{\nu}{\bar \tau^2} \phi_3) \Re\left(\overline{(\d_x+i\xi)\tau}\,(\d_x+i\xi)^2u\right)\\
&\quad
+O\big(\|u\|_{H^2_\xi}+\|\tau\|_{H^1_\xi}) (\|u\|_{H^1_\xi}+\|\tau\|_{L^2})\big)
+O(\|U\|_{H^1_\xi}\|F\|_{H^1_\xi})
\\
&\leq
-\eta\,\mathcal{E}_\xi(U) + C(\|U\|_{L^2}^2 + \|F\|_{H^1_\xi}^2).
\ea
for some positive $\eta$ and $C$ uniform with respect to $\lambda$ and $\xi$.

Combining \eqref{re} with the easy estimate 
$$
\|U\|_{L^2}^2\ =\ \frac{1}{|\lambda|}\ |\langle U, \lambda U\rangle|\leq C|\lambda|^{-1}(\|U\|_{H^1_\xi}^2 + \|F\|_{L^2}^2),
$$
for some $C$, we 
obtain
for $\Re \lambda\geq -\eta/2$ and $|\lambda|$ sufficiently large, the estimate
$$
\|U\|_{H^1_\xi}^2\ \leq\ C\|F\|_{H^1_\xi}^2
$$
for some uniform $C$. Incidentally, since $L_\xi$ has compact resolvents hence discrete spectrum composed entirely 
of eigenvalues, the above estimate also implies that such $\lambda$ do not belong to the spectrum of $L_\xi$.

More generally, by adapting the previous computations to higher-order estimates, along the lines of the method expounded in next subsection, one proves the following result required by the analysis of \cite{JZN,JNRZ}.

\begin{proposition}[Resolvent bounds]\label{spectral-damping}
For any positive integer $s$, there exist positive $\eta$, $C$ and $R$ such that if $\lambda\in\C$ is such that $|\lambda|\geq R$ and $\Re(\lambda)\geq-\eta$ then, for any $[-\pi/\Xi, \pi/\Xi)$, we have 
$$
\lambda\notin\sigma_{H^s_{per}(0,\Xi)}(L_\xi)\qquad\textrm{and}\qquad
\|(\lambda-L_\xi)^{-1}\|_{H^s_\xi\to H^s_\xi}\leq C\,.
$$
\end{proposition}

This offers a direct replacement for \cite[Appendix~B]{JZN} \emph{without assuming any condition on the background wave $\bar U$.}

\subsection{Nonlinear damping estimate}\label{s:ndamp}
%
%
%
%

The other place where high-frequency estimates play a role in the arguments of \cite{JZN,JNRZ} is in providing a nonlinear slaving bound that shows that high-regularity norms are controlled by low-regularity 
ones and enables us to close in regularity a nonlinear iteration. With the strategy implemented above we are also able to reproduce this bound {\it without assuming the slope condition \eqref{e:Lslope}.}

To be more specific let us first warn the reader that, because of the complex spatio-temporal dynamics that take 
place around periodic waves, the 
appropriate
notion of stability is neither the standard one nor the simpler orbital stability but space-modulated stability, as recalled in the next subsection. For this reason, following \cite{JZN,JNRZ}, instead of directly estimating $\tilde U-\bar U$, where $\tilde{U}=(\tilde \tau, \tilde u)$ is a solution of \eqref{eqn:co1conslaw}, we need to introduce $(V,\psi)$ such that
\begin{equation} \label{pertvar}
V(x,t)=\tilde{U}(x-\psi(x,t),t) -\bar U(x),
\end{equation}
intending to prove that $V$ and the derivatives of $\psi$ remain small provided that they are sufficiently small initially. Mark that even if $\psi$ is initially zero, as assumed in \cite{JZN}, one may not achieve the latter goal while imposing $\psi\equiv0$. In other words a modulation in space, encoded by a space-time dependent phase is in any case needed in the argument. 
See the detailed discussions in \cite{JNRZ,R}. Our new unknowns, which have to be determined together in a nonlinear way, are then $V=(\tau,u)$ and $\psi$, and a specific educated choice, that we shall not detail here, is then needed to obtain concrete equations for those. However, let us at least mention that in constructions of \cite{JZN,JNRZ} the phase shift $\psi$ is always slow so that only high regularity control on $V$ remains to be proved. This is what we provide now. 

To do so in a precise but concise way, we set $f(\tau)=(2F^2)^{-1}\tau^{-2}$, $g(\tau)=\nu\tau^{-2}$ and $h(\tau,u)=1-\tau\,u^2$. Then $\tilde U$ in \eqref{pertvar} solves \eqref{eqn:co1conslaw} provided that $V=(\tau,u)$ and $\psi$ satisfy
\be\label{modulated-swl}
\begin{array}{rcl}
(1-\psi_x)\,V_t-L\,V&=&
-\,\psi_t (\bar U+V)_x\ +\ \begin{pmatrix}0\\-\psi_x\,h(\bar U+U)\end{pmatrix}
\ +\ \begin{pmatrix}0\\\frac{-\psi_x}{1-\psi_x}\,g(\bar\tau+\tau)(\bar u+u)_x\end{pmatrix}_x\\[0.5em]
&+&
\begin{pmatrix}0\\(g(\bar\tau+\tau)-g(\bar\tau))u_x\ +\ (g(\bar\tau+\tau)-g(\bar\tau)-g'(\bar\tau)\tau)\bar u_x\end{pmatrix}_x\\[0.5em]
&-&
\begin{pmatrix}0\\f(\bar\tau+\tau)-f(\bar\tau)-f'(\bar\tau)\tau\end{pmatrix}_x
\ +\ \begin{pmatrix}0\\h(\bar U+U)-h(\bar U)-dh(\bar U)(U)\end{pmatrix}\,.
\end{array}
\ee
Defining the modified energy
\be\label{mform}
\mathcal{E}_\psi(U):=
\int (1-\psi_x)\,\Big(
\tfrac12\phi_1 \tau_x^2 +
\tfrac12\phi_2\bar \tau^3  u_x^2 +
\phi_3 \tau u_x \Big)\,,
\ee
repeating the argument of Section \ref{s:damp}, absorbing nonlinear terms into the linear ones and separating out $\psi$ terms using Sobolev's embeddings in Gagliardo-Nirenberg form and Young's inequality,
we obtain, in analogy to \eqref{dampest}, that solutions to \eqref{modulated-swl} satisfy the nonlinear estimate
\[
\frac{d}{dt}\mathcal{E}_\psi(V) \leq -\eta\,\mathcal{E}_\psi(V) 
\ +\ C\left( \|V\|_{L^2}^2+\|(\psi_t, \psi_x)\|_{H^1}^2 \right),
\]
for some positive $C$ and $\eta$, provided that we know in advance some sufficiently small upper bound on $\|(V,\psi_t, \psi_x)\|_{H^1}$ and thus are allowed to use Lipschitz bounds for $f$, $g$ and $h$ and their derivatives on a fixed neighborhood of $\bar U$. Differentiating the equations and performing the same estimate on $\partial_x^kV$, with higher-order interpolation inequalities, we obtain likewise when $k$ is a positive integer
\be\label{ndampest}
\frac{d}{dt}\mathcal{E}_\psi(\partial_x^k V) \leq -\eta\,\mathcal{E}_\psi(\partial_x^k V) 
\ +\ C\left( \|V\|_{L^2}^2+\|(\psi_t, \psi_x)\|_{H^k}^2 \right),
\ee
so long as $\|V\|_{H^1}$ and $\|(\psi_{x},\psi_t)\|_{H^k}$ remain sufficiently small.

Applying Gronwall's inequality and recalling that $\mathcal{E}_\psi(\partial_x^k V)\sim 
\|\partial_x^k V\|_{L^2}^2$ modulo lower-order terms, with constants uniform with respect to $\psi_x$ satisfying constraints above, we obtain the following key estimate showing that higher Sobolev norms $\|V\|_{H^k}$ are slaved to $\|V\|_{L^2}$ and
$\|(\psi_t,\psi_{x})\|_{H^k}$, the final nonlinear estimate needed for the analysis of \cite{JZN,JNRZ}. This provides a result analogous to \cite[Proposition~2.5]{JNRZ} and directly replacing \cite[Appendix~A]{JZN}, without making any use of a pointwise symmetrization {\it hence dropping the slope constraint \eqref{e:Lslope}.}

\begin{proposition}[Nonlinear damping]\label{damping}
For any positive integer $s$ there exist positive constants $\theta$, $C$ and $\varepsilon$ such that if $V$ and $\psi$ solve \eqref{modulated-swl} on $[0,T]$ for some $T>0$ and
$$
\sup_{t\in[0,T]}\|(V,\psi_t,\psi_x)(t)\|_{H^s(\RM)}\leq\varepsilon
$$
then, for all $0\leq t\leq T$,
\be\label{Ebds}
\|v(t)\|^2_{H^s(\RM)}
\leq Ce^{-\theta t}
\|v(0)\|_{H^s(\RM)}^2+
C\int_0^t e^{-\theta(t-s)}
\left(\|v(s)\|^2_{L^2(\RM)}+\|(\psi_t,\psi_x)(s)\|_{H^{s}(\RM)}^2\right)\dd s\,.
\ee
\end{proposition}

\subsection{Asymptotic stability}\label{s:behavior}
As discussed with great detail in \cite[Appendix~D]{JNRZ}, uniform resolvent bounds of Proposition~\ref{spectral-damping} and nonlinear slaving estimates of Proposition~\ref{damping} are the only structural conditions needed to apply almost word-by-word the arguments of \cite{JNRZ} to a periodic wave of a given 'parabolic' system. Our foregoing analysis shows that system~\eqref{swl} satisfies those around any given wave so that all conclusions of \cite{JNRZ} apply to any spectrally-stable periodic wave of \eqref{swl}. In particular, any spectrally-stable roll-wave is also nonlinearly-stable, provided that one uses definitions of stability adapted to 
periodic waves of parabolic systems, 
as we now briefly recall.

A given periodic wave solution to \eqref{swl} $\bar U$, of period $\Xi$, is said
to be
 {\it diffusively spectrally stable} provided that the generator $L$ of the linearized evolution and its Bloch symbols $L_\xi$, as defined in Subsection~\ref{s:hf}, satisfy
\begin{enumerate}
  \item[\bf{(D1)}]
$\sigma(L)\subset\{\lambda\ |\ \Re \lambda<0\}\cup\{0\}$.
  \item[\bf{(D2)}]
There exists $\theta>0$ such that for all $\xi\in[-\pi/\Xi,\pi/\Xi)$ we have
$\sigma(L_{\xi})\subset\{\lambda\ |\ \Re \lambda\leq-\theta|\xi|^2\}$.
  \item[\bf{(D3)}]
$\lambda=0$ is an eigenvalue of $L_0$ with generalized eigenspace of dimension $2$.
  \item[\bf{(H)}]
With respect to the Floquet exponent $\xi$, derivatives at $0$ of the two spectral curves passing through zero are distinct. 
\end{enumerate}
From the pioneering work \cite{S1,S2} to the recent \cite{JNRZ}, conditions (D1)--(D3) have slowly emerged as essentially sharp spectral stability conditions for periodic waves of dissipative systems. Some form of (H) is also needed but the present form could well be slightly relaxed in a near future, see precise discussion in \cite[Chapter~5]{R}. All together, conditions (D1)--(D3) and (H) express that the spectrum of $L$ is as noncritical and nondegenerate
 insofar 
as allowed by the presence around $\bar U$ of a two-dimensional family of periodic waves.

The spatial complexity of the periodic background $\bar U$ precludes any hope for a simple notion of nonlinear stability. 
Over the years there has arisen
the concrete remedy implemented in \eqref{pertvar}, 
consisting in introducing a space-time dependent phase shift, though with various possible strategies in the prescription of 
separate-- but coupled-- equations for $V$ and $\psi$. 
One obvious inspiration
for introducing a phase in the nonlinear study comes from classical analysis of simpler, asymptotically-constant patterns such as fronts, kinks, solitary waves or shock waves, for which the relevant notion of stability-- orbital stability-- already requires the introduction of a time-dependent phase. As formalized in \cite{JNRZ} the corresponding notion of stability for periodic waves-- space-modulated stability-- is obtained by measuring proximity of a function $u$ from a function $v$ in a given functional space $X$ with
$$
\delta_X(u,v)\ =\ \inf_\Psi\quad \|u\circ\Psi-v\|_X\ +\ \|\d_x(\Psi-\Id)\|_X.
$$
and not with $\|u-v\|_X$. At a given time this allows for a space-dependent phase synchronization provided that the synchronization differs from the identity by a sufficiently slow phase shift. The interested reader is again referred to \cite{JNRZ,R} for a detailed discussion of this concept. However, we stress again here that there is no hope for a better notion of stability unless the original system exhibits some nongeneric null conditions, 
denoted ``phase uncoupling'' in \cite{JNRZ}.

With these definitions in hands, our analysis combined with 
the
arguments of \cite{JNRZ} yield the following stability result.

\begin{theorem}[Nonlinear stability]\label{main-stability}
For any integer $K$, $K\geq4$, a diffusively spectrally stable periodic wave of \eqref{swl} is nonlinearly asymptotically stable from $L^1(\R)\cap H^k(\R)$ to $H^K(\R)$ in a space-modulated sense. 

More explicitly, if $\bar U$ satisfies (D1)--(D3) and (H) then, for any $K\geq4$, there exist positive $\varepsilon$ and $C$ such that any $\tilde U_0$ such that $\delta_{L^1\cap H^K}(\tilde U_0,\bar U)\leq \varepsilon$ generates a global solution $\tilde U$ to \eqref{swl} such that
$$
\forall\ t\in\R_+\,,\quad\delta_{H^K}(\tilde U(\cdot,t),\bar U)\leq\ C\ \delta_{L^1\cap H^K}(\tilde U_0,\bar U)
$$
and 
$$
\delta_{H^K}(\tilde U(\cdot,t),\bar U)\ \stackrel{t\to\infty}{\longrightarrow}\ 0\,.
$$
\end{theorem} 

The actual proof provides a much more precise statement including, for instance, a bound of $\delta_{L^p}(\tilde U(\cdot,t),\bar U)$ by $C\,(1+t)^{-\tfrac12(1-1/p)}$ similar to those for $L^p$-norms of a heat kernel or of self-similar solutions of viscous Burgers' equations. The reader is referred to \cite[Theorem~1.10]{JNRZ} for such a precise statement. Once Theorem~\ref{main-stability} is 
proved, the second part of the analysis \cite{JNRZ} may also be applied to \eqref{swl}. This yields a very precise description of the large-time asymptotic behavior in terms of a slow modulation in local parameters varying near constant parameters of the original wave, and obeying some averaged system of partial differential equations, as derived to various order of precision in \cite{W,Se,NR,JNRZ}. For 
the sake of conciseness, 
we do not state such a result here but rather refer the reader to \cite[Theorem~1.12]{JNRZ} and accompanying discussions in \cite{JNRZ,R}.

\section{The shock wave case}\label{s:shock}
As a sample of the potential wider use of the strategy expounded here, we next turn to the connection with viscous shock theory, showing
that the same linear damping estimate \eqref{dampest}
may be obtained by essentially the same argument in the asymptotically-constant, viscous shock wave case,
thus recovering the bounds established in \cite{Z1,Z2,GMWZ} by related but slightly different weighted Kawashima-type
energy estimates.\footnote{
	The weights used in \cite{Z1,Z2,GMWZ} are effectively $\phi_1=\phi_2\gg\phi_3$,
$(\phi_1)_x =-C \Big(\frac{\alpha \bar \tau^2}{\nu}
-
I\big(\frac{\alpha \bar \tau^2}{\nu} \big) \Big)\phi_1$, $C\gg1$.}
The equations of isentropic gas dynamics in Lagrangian coordinates, expressed in
a comoving frame are
\ba\label{igas}
	\partial_t \tau -c\partial_x \tau -\partial_x u&=0,
	\\
	\partial_t u -c\partial_x u +\partial_x p(\tau)
	&= \nu\partial_x(\tau^{-1}\partial_x u),
\ea
where $\tau$ is specific volume, $u$ is velocity, and $p$ is pressure.

Traveling waves $(\tau,u)(x,t)=(\bar \tau, \bar u)(x)$ satisfy the profile ODE
$$
-c^2 \tau - p(\tau) + q= cu'/\tau, \qquad q=\const.
$$
We note as in the periodic case that $c\neq 0$, else $u, p(\tau)\equiv \const$, yielding $\tau\equiv \const$, a trivial solution.
Assume that the shock is noncharacteristic, i.e., $-p'(\tau_\pm)\neq c^2$, hence $\tau_\pm$ are nondegenerate equilibria and the shock profile decays exponentially to its endstates as
$x\to \pm \infty$.

The linearized equations are
\ba \label{e:slin}
\tau_t-c\tau_x - u_x= 0,\qquad
u_t-cu_x -(\alpha \tau)_x 
=
\nu (\bar \tau^{-1}u_x)_x,
\ea
where
$\alpha:= p'(\bar \tau) + \nu \frac{\bar u_x}{\bar \tau} $.
Define $I\big(\frac{\alpha \bar \tau^2}{\nu} \big)$ to be a smooth interpolant between $\frac{\alpha \bar \tau^2}{\nu}|_{x=\pm \infty}$ such that
\be\label{expdecay}
\frac{\alpha \bar \tau^2}{\nu}
-I\big(\frac{\alpha \bar \tau^2}{\nu} \big)
=O(e^{-\theta |x|})
\ee
for some positive $\theta$.

Taking as before
$\displaystyle\mathcal{E}(U):=
\int \Big(
\tfrac12\phi_1 \tau_x^2 +
\tfrac12\phi_2\bar \tau^3  u_x^2 +
\phi_3 \tau u_x \Big)$,
we find again
\ba\label{sdiss}
\frac{d}{dt}\mathcal{E}(U(t))&=
\int\Big(
-\left(\tfrac c 2 (\phi_1)_x + \alpha \phi_3\right)\tau_x^2
-\left(\tfrac{\nu} {\bar \tau^2} \phi_2\right)u_{xx}^2
+ \left(\phi_1 -\alpha\phi_2 + \tfrac{\nu}{\bar \tau^2} \phi_3\right) \tau_x u_{xx} \Big)\\
&\quad
+ O\big((\|u\|_{H^2}+\|\tau\|_{H^1}) (\|u\|_{H^1}+\|\tau\|_{L^2})\big).
\ea
Taking
$\frac c 2 (\phi_1)_x + \Big(\frac{\alpha \bar \tau^2}{\nu}
-
I\big(\frac{\alpha \bar \tau^2}{\nu} \big) \Big)\phi_1=0$,
$\phi_1 -\alpha\phi_2 + \frac{\nu}{\bar \tau^2} \phi_3=0$,
$\phi_1(0)>0$,
and
$0<\phi_2\equiv \const \ll 1$,
we thus have
\ba\label{sgooddiss}
\frac{d}{dt}\mathcal{E}(U(t))&=
-\int
\left[ 
\left(I\big(\tfrac{\bar \tau^2 \alpha}{\nu}\big) \phi_1  -\tfrac{\alpha^2\bar \tau^2}{\nu}\phi_2\right)
\tau_x^2
+\left(\tfrac{\nu} {\bar \tau^2} \phi_2\right)u_{xx}^2\right]\\
&\quad
+ O\left((\|u\|_{H^2}+\|\tau\|_{H^1}) (\|u\|_{H^1}+\|\tau\|_{L^2})\right)\\
&\leq
-\eta'(\|u_{xx}\|_{L^2}^2 + |\tau_x\|_{L^2}^2)
+ C'(\|u\|_{H^2}+\|\tau\|_{H^1}) (\|u\|_{H^1}+\|\tau\|_{L^2}),
\ea
for some positive $C'$ and $\eta'$, and thereby the same linear damping estimate 
as in the periodic-coefficient case:
\be\label{sdampest}
\frac{d}{dt}\mathcal{E}(U(t))\leq -\eta\,\mathcal{E}(U(t)) + C\|U(t)\|_{L^2}^2
\ee
for some positive $\eta$ and $C$.

As in the periodic case, a crucial point is that 
$\displaystyle\int^x \Big(\frac{\alpha \bar \tau^2}{\nu} - I\big(\frac{\alpha \bar \tau^2}{\nu} \big) \Big)$,
hence $\phi_1$ and $1/\phi_1$, remains bounded, so that
$\mathcal{E}(U)\sim (\|u_{x}\|_{L^2}^2 + \|\tau_x\|_{L^2}^2)$ modulo $\|\tau\|_{L^2}^2$,
a property following in this case by exponential decay, \eqref{expdecay}.


The above may be recognized as exponentially weighted Kawashima-type estimates similar to those used in the study of viscous shock stability in \cite{Z1,Z2,GMWZ}, reflecting the growing analogy between the periodic and asymptotically-constant cases. Actually, part of the recent activity of the authors, jointly with others, focused on dynamics around periodic waves, and 
culminating more or less
in \cite{JNRZ}, was motivated by the will to put its analysis on a par with classical ones on asymptotically-constant waves. With this respect the present contribution that provides analytical tools necessary to consider large-amplitude periodic waves of hyperbolic-parabolic systems should be compared with \cite{MaZ}, where, motivated by some clever "transverse" energy estimates of Goodman \cite{G} in the study of small-amplitude stability, the treatment of large-amplitude viscous shock waves was first carried out. In the reverse direction the Conjugation Lemma of \cite{MZ} on asymptotically constant-coefficient coordinate transformations from asymptotically constant- to constant-coefficient systems may be thought as analogous to the classical Floquet Lemma on periodic coordinate transformation of periodic- to constant-coefficient systems of equations. 

\section{Discussion}\label{s:discussion}

At the linear level, for a general second-order hyperbolic-parabolic principal part $U_t + AU_x=(BU_x)_x$, $\Re B\geq 0$, a Kawashima-type estimate is on an energy combining $\mathcal{E}(U):=\langle U_x,A^0 U_x\rangle + \langle U,K U_x\rangle$,
with the lower-order $\langle U,A^0 U\rangle$, where $A^0$ is symmetric positive definite and $K$ is skew symmetric, chosen, where possible, so that 
\be\label{co}
\Re(A^0B+KA)>0.
\ee
When $A$ and $B$ are constant, as arising from linearization around a constant state, and the original nonlinear system admits a strictly convex entropy, \eqref{co} may be reduced to a simple-looking condition that is satisfied by most of systems of physical interest; see \cite{Kawashima-thesis,Liu_Zeng}. As a result, for small-amplitude waves, a suitable choice of $K$ may typically be achieved globally with a constant $K$. However, for large-amplitude shocks, this can be done typically only near $x\to \pm \infty$ where $A$ is symmetrizable, and one needs to recover coercivity in the near field $|x|\leq C$ in a different way.

As exemplified here-- for the first time in a periodic context, the key to the treatment of large-amplitude waves, is to choose the "symmetrizer" $A_0$ jointly with the "compensator" $K$ so that one may use a clever choice for $A_0$ to relax constraints on $K$ and \emph{vice versa}. In the present case, 
$$
A=\bp -c & -1\\ \alpha & -c\ep,
$$
for either of the St. Venant or isentropic compressible Navier--Stokes equations. The issue in the latter case is that symmetrizability holds in general only in the limits $x\to \pm \infty$, in the former that it holds only on average, but in any case not pointwise.
However, we have seen that energy estimates can be recovered by modulating classical symmetrizers and compensators with appropriate asymptotically-constant, or periodic exponential weights.

Mark that our analysis, while apparently quite robust, leaves widely open the question of determining, for general systems, what kind of notion of symmetrizability on average could lead to similar periodic-coefficient high-frequency damping.

\bibliographystyle{alpha}
\bibliography{Ref_slope}

\end{document}